\documentclass[11pt,twoside]{article}

\usepackage{latexsym, amssymb, amsmath, theorem}

\numberwithin{equation}{section}

\theoremstyle{change}
\theorembodyfont{\itshape}
\newtheorem{theorem}{Theorem}[section]
\newtheorem{proposition}[theorem]{Proposition}
\newtheorem{lemma}[theorem]{Lemma}

\theorembodyfont{\rmfamily}
\newtheorem{definition}[theorem]{Definition}

\newtheorem{example}[theorem]{Example}
\newtheorem{remark}[theorem]{Remark}
\newenvironment{proof}{{\noindent \textbf{Proof}\,\,}}{\hspace*{\fill}$\Box$\medskip}

\begin{document}
{\bf Syst{\`e}mes dynamiques \slash Dynamical systems}

\

\

\centerline{{\bf \Large{On convergence of generalized continued
      fractions}}} 
\centerline{{\bf \Large{and Ramanujan's conjecture}}}

\

\centerline{\large{A.A.Glutsyuk}}

\def\cc{\mathbb C}
\def\var{\varepsilon}

\begin{abstract} We consider continued fractions 
\begin{equation} 
\frac{-a_1}{1-\frac{a_2}{1-\frac{a_3}{1-\dots}}}
\label{fr}\end{equation}
with real coefficients $a_i$ converging to a limit $a$.
S.Ramanujan had stated the theorem (see [ABJL], p.38) 
saying that if $a\neq\frac14$, then the fraction converges if and only if 
$a<\frac14$. The statement of convergence was proved in 
[V] for complex $a_i$ converging to $a\in\mathbb C\setminus[\frac14,+\infty)$ 
(see also [P]). J.Gill [G] proved the divergence of (\ref{fr}) under the assumption 
that $a_i\to a>\frac14$ fast enough, more precisely, whenever 
\begin{equation}\sum_i|a_i-a|<\infty.\label{gill}\end{equation}
The Ramanujan's conjecture saying that (\ref{fr}) diverges always whenever 
$a_i\to a>\frac14$ remained up to now an open question. 
In the present paper we disprove it. We show (Theorem \ref{th1}) that for any 
$a>\frac14$ there exists a real sequence $a_i\to a$ such that (\ref{fr}) 
converges\footnote{The 
author acknowledges that Alexey Tsygvintsev has constructed 
(by a completely different method) 
a beautiful explicit example 
[Ts] of a sequence $a_i\to1$ given by a simple recurrent formula. This example 
comes from the analytic function theory}. Moreover, we show (Theorem \ref{go}) 
that Gill's sufficient divergence condition (\ref{gill}) is the optimal 
condition on the speed of convergence of the $a_i$'s. 
\end{abstract}

\centerline{{\bf \large{Sur convergence des fractions continues 
g{\'e}n{\'e}ralis{\'e}es}}} 
\centerline{{\bf \large{et  
une conjecture de Ramanujan}}}

\

{\bf R{\'e}sum{\'e}.} Nous consid{\'e}rons une fraction continue 
$$ 
\frac{-a_1}{1-\frac{a_2}{1-\frac{a_3}{1-\dots}}}
\label{frf} \ \ \ \ \ (0.1)$$
{\`a} co{\'e}fficients r{\'e}els $a_i\to a$. S.Ramanujan a formul{\'e} 
le th{\'e}or{\`e}me (voir [ABJL], p.38), qui disait, que si $a\neq\frac14$, 
alors la fraction converge,  si et seulement si $a<\frac14$. La 
convergence a {\'e}t{\'e} d{\'e}montr{\'e}e dans [V] pour $a_i$ complexes 
convergeant vers un $a\in\cc\setminus[\frac14,+\infty)$ (voir aussi 
[P]). J.Gill [G] a d{\'e}montr{\'e}, que la fraction diverge, si 
$a_i\to a>\frac14$ assez vite, plus pr{\'e}cisement, si 
$\sum_i|a_i-a|<\infty$. 

La conjecture de Ramanujan disant que 
la fraction diverge toujours, quand $a>\frac14$, restait ouverte 
jusqu'au pr{\'e}sent.  Nous 
montrons, qu'elle est fausse: pour tout $a>\frac14$ il existe une suite 
r{\'e}elle $a_i\to a$ telle que la fraction converge (Th{\'e}or{\`e}me 1.1).  
Nous montrons aussi, que la condition pr{\'e}cedante de Gill, 
qui est suffisante pour que la fraction diverge, est celle 
optimale sur la vitesse de convergence des $a_i$. 

\centerline{{\bf \Large{Version fran\c caise abr{\'e}g{\'e}e}}}

\

{\bf Th{\'e}or{\`e}me 1.1.} Pour tout $a>\frac14$ il existe une suite 
$a_i\to a$, telle que la fraction (0.1)  converge. 

{\bf Th{\'e}or{\`e}me 1.2.} {\'E}tant donn{\'e}s un $q\in\mathbb N$, $q\geq3$, 
et une suite $r_i\to0$, $r_i>0$, telle que 
$\sum_ir_i=\infty$. Alors il existe toujours un $a>\frac14$ et une 
suite  
$$a_i\to a,\ a_i=a,\ \text{si} 
\ i\not\equiv1,2(modq),\ a_{qi+1}-a=O(r_i), a_{qi+2}-a=O(r_i),\ 
\text{quand}\ i\to\infty,$$
 telle que la fraction (0.1) converge.

\section{Main results and the plan of the paper}

\subsection{Main results}
\begin{theorem} \label{th1} For any $a>\frac14$ there exists a 
real sequence 
$a_i\to a$ such that the continued fraction (\ref{fr}) converges. 
\end{theorem}

\begin{theorem} \label{go} Given any $q\in\mathbb N$, $q\geq3$, 
and a sequence $r_i\to0$, $r_i>0$, such that 
\begin{equation}\sum_i{r_i}=\infty\label{rash}.\end{equation}
Then there exists an $a>\frac14$ and a real sequence 
\begin{equation}a_i\to a,\ a_i=a,\ \text{if} 
\ i\not\equiv1,2(modq),\ a_{qi+1}-a=O(r_i), a_{qi+2}-a=O(r_i),\ 
\text{as}\ i\to\infty,\label{fgo}\end{equation}
such that the continued fraction (\ref{fr}) converges. 
\end{theorem}

\subsection{The plan of the proofs and generalizations}
\def\rr{\mathbb R}
\def\nn{\mathbb N}
\def\zz{\mathbb Z}

In the proof of Theorem \ref{th1} we use the following expression of 
continued fraction (\ref{fr}) as a limit of compositions of 
M{\"o}bius transformations of the closed upper half-plane 
$H=\{ Im z\geq0\}$. For any $b\in\rr$ 
define the M{\"o}bius transformation 
$$T_b:H\to H,\ T_b(z)=-\frac b{z+1}.$$
\begin{proposition} The subsequent ratios $\frac{p_n}{q_n}$ of 
the continued fraction (\ref{fr}) are given by the formula 
\begin{equation}\frac{p_n}{q_n}=\tau_n=T_{a_1}\circ\dots\circ T_{a_n}(0).
\label{tn}\end{equation}
\end{proposition}
The Proposition is well-known and follows immediately from 
definition (by induction in $n$). 

Recall the following 

\begin{definition} A M{\"o}bius transformation $M:H\to H$ of the 
upper half-plane is said to be {\it elliptic} (respectively, hyperbolic), 
if it has a fixed point in $Int(H)$ (respectively, two fixed points 
on the boundary of $H$: then one of them is attractor, the other one is 
repeller). (Each elliptic element is M{\"o}bius conjugated to a rotation  
of unit disc. By definition, its {\it rotation number} is 
$(2\pi)^{-1}$ times the corresponding 
rotation angle.) The {\it multiplier} of a hyperbolic transformation 
$T$ (denoted $\mu=\mu(T)$) is its 
derivative at the attractor (by definition, $0<\mu<1$). 
\end{definition}

\def\qq{\mathbb Q}

\begin{remark} \label{c'a} The transformation $T_a$ is elliptic if and only if $a>\frac14$. 
The function $a\mapsto\rho(a)$, whose value is the rotation number of $T_a$, 
is an analytic diffeomorphism $(\frac14,+\infty)\to(0,\frac12)$. 
The function $a\mapsto c(a)$ whose value is the fixed point of $T_a$ in $H$ is 
an analytic diffeomorphism $(\frac14,+\infty)\to i\rr_+-\frac12$. One has 
$$\rho(1)=\frac13,$$
since $T_1$ permutes cyclically 0, 1 and $\infty$. If $\rho(a)=\frac pq\in\qq$, 
then $T_a^q=Id$. 
\end{remark}

\begin{example} Let $a_i\equiv a>\frac14$. Then $\tau_n=T_a^n(0)$ 
does not have a limit: this is either periodic or a quasiperiodic sequence.
\end{example}

\begin{theorem} \label{th2} Theorem \ref{th1} holds for 
each $a$ with $\rho(a)\in\qq$. 
Moreover, given any smooth  
elliptic transformation family $\{T_a\}_{a\in\rr}: H\to 
H$ with a smooth fixed point family $c(a)$. Let $a$ 
be a parameter value such that $\rho(a)=\frac 
pq\in\qq\setminus\{0,\frac12\}$ 
and $c'(a)\neq0$. Then 
there exists a sequence $a_i\to a$ such that the corresponding sequence 
$\tau_n$ from (\ref{tn}) converges. For any given diverging series 
$\sum_ir_i=\infty$, $r_i\to0_+$, the sequence $a_i$ can be chosen to 
satisfy (\ref{fgo}).  
\end{theorem}

\begin{theorem} \label{th3} Theorem \ref{th1} holds for each $a$ with $\rho(a)
\notin\qq$. Moreover, let $T_a$ be an elliptic transformation family as in  
Theorem \ref{th2}. Let $a$ be a parameter value such that 
$\rho(a)\notin\qq$ and $\rho\not\equiv const$ near $a$.  Then there 
exists a sequence $a_i\to a$ such that $\tau_n$ converge.
 \end{theorem}

Theorems \ref{th2} and \ref{th3} imply Theorem \ref{th1}. Theorem 
\ref{th2} implies Theorem \ref{go}. It is proved 
in the next Section. Theorem \ref{th3} is proved in Section 3.

\def\a{\alpha}
\def\b{\beta}

\section{Limits with rational $\rho(a)$. Proof of Theorem \ref{th2}}

Let $\rho(a)=\frac pq\in\qq$, thus, $T_a^q=1$. 
We choose appropriate sequences $\a_r,\b_r\to a$ as specified below and 
put  
\begin{equation}a_i=a\ \text{if} \ i\not\equiv 1,2 (mod q); \ 
a_{qr+1}=\a_r,\ a_{qr+2}=\b_r.\ \text{Denote}\label{defai}\end{equation} 
$$T_{\a,\b,q}=T_{\a}\circ T_{\b}\circ T_a^{q-2},\ \a,\b\in\rr. 
\ \text{We choose}\ \a_r,\b_r\ \text{so that}$$ 
\begin{equation}T_{\a_r,\b_r,q} 
 \ \text{be hyperbolic, denote}\ A_r,\ R_r 
\ \text{their attractors (respectively, repellers)},  
\label{hyp}\end{equation}
\begin{equation}A_i\to A, R_i\to R, 
\ A\neq R,\ R\notin C_a=\{ T_a^l(0), l=0,\dots,q-1\},
\label{ar}\end{equation}
\begin{equation}\prod_r\mu_r=0,\ \mu_r=\mu(T_{\a_r,\b_r,q}).
\label{mu}\end{equation}

The possibility of the above choice of $\a_r$, $\b_r$ is proved at the end of 
the Section.

Below we show (the next Proposition and the paragraph after) that 
the corresponding sequence $\tau_n$ converges, whenever 
conditions (\ref{hyp})-(\ref{mu}) hold. 

\begin{proposition} Let $H_1$, $H_2,\dots$ be an arbitrary sequence of hyperbolic 
transformations $H\to H$ of the upper half-plane, $A_i$, $R_i$ be respectively 
their attractors and repellers, $A_i\to A$, $R_i\to R\neq A$. Let 
(\ref{mu}) hold with $\mu_r=\mu(H_r)$. Then the mapping sequence 
$\widehat H_n=H_1\circ\dots\circ H_n$ converges uniformly to a constant 
mapping on compact sets in $\partial H\setminus R$. 
\end{proposition}
\begin{proof} If $A_i=A$, $R_i=R$, then the statement of the Proposition 
follows immediately. If we fix a compact set $K\Subset\partial H\setminus 
\{ A\cup R\}$, then for any $i$ large enough the transformation 
$H_i$ moves the points of $K$ towards $A$ (along $\partial H$) by  
asymptotically the same distance, as the hyperbolic transformation with the 
same multiplier $\mu_i$ but with $A_i=A$, $R_i=R$. This together with 
the previous statement and the monotonicity of the restrictions $H_i|_{\partial H}$ 
implies the Proposition. 
\end{proof}
\subsection{Proof of convergence of $\tau_n$} 
If (\ref{hyp})-(\ref{mu}) hold, then the transformations 
$H_r=T_{\a_r,\b_r,q}$ satisfy the conditions of the previous Proposition. 
Hence, their compositions 
$\widehat H_r=T_{a_1}\circ T_{a_2}\circ\dots\circ T_{a_{qr}}$ converge uniformly 
to a constant limit (denote it $x$) on compact sets 
in $\partial H\setminus R$. By definition, $\tau_{qr}=\widehat H_r(0)$, 
$0\neq R$ by (\ref{ar}). This implies that $\tau_{qr}\to x$, as $r\to\infty$.
To show that the whole sequence $\tau_n$ converges to $x$, 
we use condition (\ref{ar}), which says that the finite $T_a$- orbit 
$C_a$ of 0 does not meet $R$. Let us take a $\delta>0$ so that the closed 
$\delta$- neighborhood $U$ of the latter orbit be disjoint from $R$. Then 
$\widehat H_r\to x$ uniformly on $U$. If 
$r$ is large enough, then  
$$T_{\a_r}(0), T_{\a_r}\circ T_{\b_r}(0), T_{\a_r}\circ T_{\b_r}\circ 
T_a^{l-2}(0)\in U 
\ (0<l\leq q). \ \text{By definition,}$$ 
 $$\tau_{qr+1}=\widehat H_r\circ T_{\a_r}(0), 
\tau_{qr+l}=\widehat H_r\circ T_{\a_r}\circ T_{\b_r}\circ T_a^{l-2}(0).$$ 
The two last statements imply together that the $q$ sequences 
$\tau_{qr},\tau_{qr+1},\tau_{qr+2},\dots,\tau_{qr+q-1}$ converge to $x$, hence, 
the whole sequence $\tau_n$ converges. Theorem \ref{th2} is proved. 
\subsection{The construction of sequences $\a_i$, $\b_i$ satisfying  (\ref{hyp})-(\ref{mu})} 

\begin{lemma} \label{mtl}{\bf (Main Technical Lemma).} Let a family 
$T_a$ and a parameter value $a$ be as 
in Theorem \ref{th2}.  Then for any point $R\in\partial H$ 
(maybe except two points) there exist two linear families of parameter values 
\begin{equation}\a(t)=a+c_1t,\ \b(t)=a+c_2t,\ c_1,c_2\in\rr,\label{abt}\end{equation} 
such that for any $t>0$ small enough the transformation 
$T_{\a(t),\b(t),q}=T_{\a(t)}\circ T_{\b(t)}\circ T_a^{q-2}$ be hyperbolic 
and its repeller $R(t)$ (respectively, 
attractor $A(t)$) tends to $R$ (respectively, to a point $A\neq R$), as 
$t\to 0_+$. 
Moreover, one can achieve that the families $A(t)$, $R(t)$ be smooth at 0, 
and the derivative in $t$ at $t=0$ of the previous family $T_{\a(t),\b(t),q}$ 
be nonzero. 
\end{lemma}

The Lemma is proved below. 

Let $\rho(a)=\frac pq$, $C_a$ be the (finite) $T_a$- orbit of $0$. Let us
choose any $R\notin C_a$ that satisfies the statements of Lemma \ref{mtl}.
Let $\a(t)$, $\b(t)$ be the corresponding families (\ref{abt}).  Take a
sequence $t_k\to0_+$ and put $\a_k=\a(t_k)$, $\b_k=\b(t_k)$.  The
conditions (\ref{hyp}) and (\ref{ar}) follow immediately from
construction. Condition (\ref{mu}) holds, if and only if $\sum t_i=\infty$
(these are the $t_i$ we choose). This follows from the fact that the
function $\mu(t)=\mu(T_{\a(t),\b(t),q})$ has nonzero derivative at 0:
$\mu(t)=1+st+O(t^2)$, $s\neq0$, hence, $\ln\mu_k=-st_k(1+o(1))$.  Indeed,
otherwise, if $\mu'(0)=0$, the transformation family $T_{\a(t),\b(t),q}$
would have zero derivative in $t$ at $t=0$ - a contradiction to the last
statement of the Lemma. This finishes the construction. Statement
(\ref{fgo}) follows immediately, if we put $t_i=r_i$. Theorem \ref{th2} is
proved.

\begin{proof} {\bf of Lemma \ref{mtl}.} Consider $T_{\a,\b, q}=T_{\a}\circ T_{\b}\circ T_a^{q-2}$ 
as a family of mappings depending on two variable parameters $\a$ and $\b$ 
(the $a$ is fixed). It is identity, if $\a,\b=a$. Consider its derivative 
in $\a$ at $(\a,\b)=(a,a)$ (which is a vector field on $\partial H$ denoted $v_1$). 
Its derivative in $\b$ at the same point 
is another vector field $v_2=(T_a)_{*}v_1$ that is 
the image of $v_1$ under the diffeomorphism $T_a:\partial H\to\partial H$. 
We claim that the fields $v_1$ and $v_2$ are not constant-proportional. 
Indeed, otherwise the group generated by $T_a$ and the 1-parametric 
subgroup in $Aut(H)$ generated by $v_1$ 
would be solvable. Since $T_a$ is elliptic, this 
implies that $T_a$ either belongs to the same 
1-parameter subgroup, or is an involution. The first case 
is impossible: otherwise the centers $c(a)$ would have zero derivatives at $a$ - 
a contradiction to the conditions of the Lemma. The second case is impossible by 
the hypothesis $\rho(a)\neq0,\frac12$. 
 
Thus, the vector fields $v_1$ and  $v_2$ are not proportional. Hence, 
for any point $R\in\partial H$ one can find 
a linear combination $v=c_1v_1+c_1v_2\not\equiv0$ that vanishes at $R$. If the 1- jet  
of $v$ at $R$ does not vanish (then one can achieve that 
$v'(R)>0$ by changing sign), this implies that 
$v$ has another zero $A\in\partial H\setminus R$. Then the corresponding 
families (\ref{abt}) are those we are looking for. If the latter 1- jet  
vanishes, this implies that the commutator $[v_1,v_2]$ 
(which also belongs to the Lie algebra of the group $Aut(H)$) 
vanishes at $R$. The latter commutator does not vanish identically 
(since $v_i$ are not proportional) and cannot have more than two zeros. 
This together with the previous discussion proves the Lemma 
\end{proof}

\section{Case of irrational limit rotation. Proof of Theorem \ref{th3}.} 

\def\wt#1{\widetilde#1}

Let $\ \rho(a)\notin\qq,\ \wt a_n\to a,\ \rho(\wt a_n)=
\frac{p_n}{q_n}\in\qq.$ 
 We choose appropriate $\a_n,\b_n\to a$, a natural number sequence 
 $N_1,N_2,\dots$, and define $a_n$ as follows: 
 
 1) Let $n\leq N_1q_1$. Put 
  $a_n=\wt a_1,\ \text{if}\ n\not\equiv1,2(mod q_1); \ 
a_{q_1r+1}=\a_1, \  a_{q_1r+2}=\b_1,\ r<N_1$.  

2) Let $N_1q_1<n\leq N_1q_1+N_2q_2$. Put 
$n_1=n-N_1q_1,\ a_n=\wt a_2,\ \text{if} \ n_1\not\equiv1,2(mod q_2)$;  
$a_n=\a_2,\ \text{if} \ n_1\equiv1(mod q_2);\  
a_n=\b_2,\ \text{if}\ n_1\equiv2(mod q_2),\ \text{etc.}$
 
We show that (\ref{fr}) converges if we take $\a_i$, $\b_i$ and $N_i$ as 
 specified below. 
 
 {\bf Choice of $\a_n$ and $\b_n$.} Denote 
 $$\psi_0=\psi_{n,0}=Id, \ \psi_n=T_{\a_n,\b_n,q_n}=T_{\a_n}\circ 
T_{\b_n}\circ 
T_{\wt a_n}^{q_n-2}, \psi_{n,1}=T_{\a_n},\ 
\psi_{n,2}=\psi_{n,1}\circ T_{\b_n},$$ 
$\psi_{n,l}=\psi_{n,2}\circ
 T_{\wt a_n}^{l-2},\ \text{for}\ 2\leq l\leq q_n-1.$ 
  We choose $\a_n$ and $\b_n$ so that  the transformations 
\begin{equation}\psi_n\ \text{be hyperbolic
, denote}\ A_n,\ R_n\ \text{their attractors (respectively, repellers),}
\label{hypn}\end{equation} 
\begin{equation}R_n\notin M_n=\{A_{n+1},\psi_{n+1,l}(A_{n+1}),0, 
\psi_{n+1}^r\circ\psi_{n+1,l}(0)|\ 0\leq l\leq q_{n+1}-1,\ r\in\mathbb N\cup 0.\}
\label{rnmn}\end{equation}

\begin{remark} \label{remrm}
The previous set $M_n$ is infinite and accumulates exactly to 
the finite $T_{\wt a_n}$- orbit of $A_{n+1}$, which follows from definition. 
This implies that if (\ref{rnmn}) holds, then $M_n$ does not accumulate 
to $R_n$. Thus, in this case 
choosing appropriate power $N_n$, one can achieve that 
the image $\psi_n^{N_n}(M_n)$ be arbitrarily close to $A_n$. 
\end{remark}

\def\tk{\theta_k}

{\bf Choice of $N_i$.} Let $\a_i$, $\b_i$ be already chosen to satisfy 
(\ref{hypn}) and (\ref{rnmn}). Denote 
$$\theta_k=\psi_1^{N_1}\circ\dots\circ\psi_k^{N_k}.$$
We construct $N_i$ (by induction in $i$) 
in such a way that 
\begin{equation}diam(\theta_k(M_k))<\frac1{2^k}.\label{diamk}\end{equation}
The possibility to do this follows immediately from the last statement of 
the previous Remark. 
Let us show that then the sequence $\tau_n$ is Cauchy 
(hence, converges). Denote $n_k=\sum_{i=1}^kq_iN_i$. It suffices to show that
\begin{equation}\text{for any}\ k\ \text{and any}\ m\geq n_k\ \text{one has}\ 
dist(\tau_{n_k},\tau_m)<\frac1{2^{k-2}}.\label{distk}\end{equation}
Case 1: $m=n_i>n_k$, say, $m=n_{k+1}$. 
Then $\tau_{n_k}=\tk(0)$, $\tau_m=\theta_{k+1}(0)=
\tk\circ\psi_{k+1}^{N_{k+1}}(0)$. By definition, 
$0, \psi_{k+1}^{N_{k+1}}(0)\in M_k$. By (\ref{diamk}), 
$dist(\tk(0),\theta_{k+1}(0))=dist(\tk(0),\tk(\psi_{k+1}^{N_{k+1}}(0)))
<\frac1{2^k}.$
\begin{equation}
\text{Therefore,}\ 
dist(\tk(0),\theta_s(0))<\frac1{2^{k-1}}\ \text{for any}\ 
s>k.\label{dists}\end{equation}
This proves (\ref{distk}) for any $m=n_i>n_k$. 

General case: $m>n_k$ is arbitrary. Take $s\in\mathbb N$ such that 
$$n_s\leq m<n_{s+1}.\ \text{Then}\ m=n_s+rq_{s+1}+l,\ 0\leq l<q_{s+1},\ 
\tau_m=\theta_s\circ\psi_{s+1}^{r-1}\circ\psi_{s+1,l}(0).$$
Analogously to the previous discussion, by (\ref{diamk}), 
$dist(\theta_s(0),\tau_m)<\frac1{2^s}$. This together with  
 (\ref{dists}) implies (\ref{distk}). Theorem \ref{th3} is proved. 

\section{Acknowledgements} 
The author is grateful to A.Tsygvintsev for attracting his attention 
to the problem and for helpful discussions. 

\section{Bibliography}

[ABJL] G.E.Andrews, B.C.Berndt, L.Jacobsen, R.L.Lamphere. 
- The continued fractions found in the unorganized portions of Ramanujan's 
notebooks. - Memoirs of the AMS, vol. 99 (1992), No 477. 

[G] J.Gill, Infinite composition of M{\"o}bius transformations, - Trans. Amer. 
Math. Society, 176 (1973), 479-487. 

[P] O.Perron, Die Lehre von den Kettenbr{\"u}chen, - Band 2, dritte Auf., 
B.G. Teubner, Stuttgart, 1957. 

[Ts] A.Tsygvintsev. On the convergence of continued fractions at 
Runckel's points and the Ramanujan conjecture. - Preprint 
http://arxiv.org/abs/math.FA/0412298

[V] E.B.Van Vleck, On the convergence of algebraic continued fractions whose 
coefficients have limiting values. - Trans. Amer. Math. Soc. 5 (1904), 
253-262. 

\

CNRS, Unit{\'e} de Math{\'e}matiques Pures et Appliqu{\'e}es, M.R., 
{\'E}cole Normale Sup{\'e}rieure de Lyon, 46 all{\'e}e d'Italie, 
69364 Lyon Cedex 07, France

{\it E-mail address:} aglutsyu$@$umpa.ens-lyon.fr

\end{document}